\documentclass{article}
\usepackage{amssymb}
\usepackage{amsfonts}
\usepackage{amsmath}

\newtheorem{theorem}{Theorem}
\newtheorem{acknowledgement}[theorem]{Acknowledgement}

\newenvironment{proof}[1][Proof]{\noindent\textbf{#1.} }{\ \rule{0.5em}{0.5em}}
\input{tcilatex}

\begin{document}

\title{Ricci identities of the Liouville d-vector fields $z^{\left( 1\right)
\alpha }$ and $z^{\left( 2\right) \alpha }$}
\author{Oana Alexandru \\
"Transilvania" University of Brasov}
\date{}
\maketitle

\begin{abstract}
It is the purpose of the present paper to outline an introduction in theory
of embeddings in the manifold $Osc^{2}M$. First, we recall the notion of
2-osculator bundle (\cite{at},\cite{at2}). The second section is dedicated
to the notion of submanifold in the total space of the 2-osculator bundle,
the manifold $Osc^{2}M$. A moving frame is constructed. The induced N-linear
connections and the relative covariant derivatives are discussed in third
and fourth sections. The Ricci identities of the Liouville d-vector fields
are present in the last section.


\end{abstract}

\begin{center}
\bigskip \bigskip {\LARGE Introduction}
\end{center}

\bigskip Generally, the geometries of higher order defined as the study of
the category of bundles of jet $\left( J_{0}^{k}M,\pi ^{k},M\right) $ is
based on a direct approach of the properties of objects and morphisms in
this category, without local coordinates.

But, many mathematical models from Lagrangian Mechanics, Theoretical Physics
and Variational Calculus used multivariate Lagrangians of higher order
acceleration.

From here one can see the reason of construction of the geometry of the
total space of the bundle of higher accelerations (or the osculator bundle
of higher order) in local coordinates.

As far we know the theory of Finsler or Lagrange submanifolds is far from
being settled. In \cite{mian1} and \cite{mian2} R. Miron and M. Anastasiei
give the theory of subspaces in the Lagrange or generalized Lagrange spaces.
Also, in \cite{mi} R. Miron presented the theory of subspaces in higher
order Lagrange spaces.

This article is draw upon the original construction of the higher order
geo\-me\-try given by R. Miron and Gh. Atanasiu (\cite{at},\cite{mi},\cite%
{miat1996},\cite{miat},\cite{miat1994}).

In \cite{mi}, R. Miron construct the theory of the subspaces in higher order
Lagrange spaces using the\ canonical metrical N-connection of the space $%
L^{\left( 2\right) n}$ having \textbf{three} coeficients $\left( F_{jk}^{i},%
\underset{1}{C}\overset{}{_{jk}^{i}},\underset{2}{C}\overset{}{_{jk}^{i}}%
\right) .$ In our work, we take the canonical metrical N-linear connection
of the manifold $Osc^{2}M$ having \textbf{nine} coeficients $\left( \underset%
{\left( i0\right) }{L}\overset{}{_{bc}^{a}},\underset{\left( i1\right) }{C}%
\overset{}{_{bc}^{a}},\underset{\left( i2\right) }{C}\overset{}{_{bc}^{a}}%
\right) ,\left( i=0,1,2\right) .$

If $\check{M}$ is an immersed manifold in the manifold $M$, a nonlinear
connection N on the manifold $Osc^{2}M$ induce a nonlinear connection $%
\check{N}$ on\ the submanifold $Osc^{2}\check{M}.$ We take the canonical
metrical N-linear connection $D$ on the manifold $Osc^{2}M$ and we obtain
the induced tangent and normal connections. Also, we introduce the relative
covariant derivatives in the algebra of d-tensor fields (\cite{Oa1}). Next,
we get the Ricci identities for the Liouville d-vector fields $z^{\left(
1\right) \alpha }$ and $z^{\left( 2\right) \alpha }$ (Theorem 5.2). The same
problem was solved by prof. Atanasiu Gh. in \cite{at} for the Liouville
d-vector fields $z^{\left( 1\right) a}$ and $z^{\left( 2\right) a}$ on the
manifold $Osc^{2}M.$

\section{The 2-osculator bundle $\left( Osc^{2}M,\protect\pi ^{2},M\right) $}

Let M be a real differentiable manifold of dimension n, whose coordinates
are $\left( x_{a}\right) _{a=\overline{1,n}}$. Note that, throughout this
paper the indices $a,b,...$ run over set $\{1,2,...,n\}.$ The Einstein
convention of summarizing is adopted all over this work.

Let us consider two curves $\rho ,\sigma :I\rightarrow M$, $\ $having images
in a domain of local chart $U\subset M.$ We say that $\rho $ and $\sigma $
have a "contact of order 2" in a point $x_{0}\in U$ if $\rho \left( 0\right)
=\sigma \left( 0\right) =x_{0},$ $\left( 0\in I\right) ,$ and for any
function $f\in \mathcal{F}\left( U\right) $%
\begin{equation}
\frac{d^{\beta }}{dt^{\beta }}\left( f\circ \rho \right) \left( t\right)
\mid _{t=0}=\frac{d^{\beta }}{dt^{\beta }}\left( f\circ \sigma \right)
\left( t\right) \mid _{t=0},\left( \beta =1,2\right) .  \tag{1.1}
\label{1.11}
\end{equation}

The relation "contact of order 2" is an equivalence on the set of smooth
curves in M, which pas through the point x$_{0}.$ Let $[\rho ]_{x_{0}}$ be a
class of equivalence. It will be called a "2-osculator space" in a point $%
x_{0}\in M.$ The set of 2-osculator spaces in the point $x_{0}\in M$ will be
denoted by $Osc_{x_{0}}^{2}M$, and we put%
\begin{equation*}
Osc^{2}M=\underset{x_{0}\in M}{\cup }Osc_{x_{0}}^{2}M.
\end{equation*}%
One considers the mapping $\pi ^{2}:Osc^{2}M\rightarrow M$ define by $\pi
^{2}\left( \left[ \rho \right] _{x_{0}}\right) =x_{0}.$ Obviously, $\pi ^{2}$
is a surjection.

The set $Osc^{2}M$ is endowed with a natural differentiable structure,
induced by that of the manifold M, so that $\pi ^{2}$ is a differentiable
maping. It will be descrieb bellow.

The curve $\rho :I\rightarrow M$ $\left( \func{Im}\rho \subset U\right) $ is
analytically represented in the local chart $\left( U,\varphi \right) $ by $%
x_{0}=x_{0}^{a}\left( =x^{a}\left( 0\right) \right) .$ Taking the function $%
f $ from $\left( \text{\ref{1.11}}\right) $, succesively equal to the
coordinate functions $x^{a},$ then a representative of the class $[\rho
]_{x_{0}}$ is given by%
\begin{equation*}
x^{\ast a}\left( t\right) =x^{a}\left( 0\right) +t\dfrac{dx^{a}}{dt}\left(
0\right) +\frac{1}{2}t^{2}\dfrac{d^{2}x^{a}}{dt^{2}}\left( 0\right) ,\text{ }%
t\in \left( -\varepsilon ,\varepsilon \right) \subset I.
\end{equation*}%
The previous polynomials are determined by the cofficients%
\begin{equation}
x_{0}^{a}=x^{a}\left( 0\right) ,y^{\left( 1\right) a}=\dfrac{dx^{a}}{dt}%
\left( 0\right) ,y^{\left( 2\right) a}=\frac{1}{2}\dfrac{d^{2}x^{a}}{dt^{2}}%
\left( 0\right) .  \tag{1.2}  \label{1.22}
\end{equation}%
Hence, the pair $\left( \left( \pi ^{2}\right) ^{-1}\left( U\right) ,\Phi
\right) ,$ with $\Phi \left( \lbrack \rho ]_{x_{0}}\right) =\left(
x_{0}^{a},y^{\left( 1\right) a},y^{\left( 2\right) a}\right) \in R^{3n},$ $%
\forall \lbrack \rho ]_{x_{0}}\in \left( \pi ^{2}\right) ^{-1}\left(
U\right) $ is a local chart on $Osc^{2}M$. Thus a differentiable atlas $%
\mathcal{A}_{M}$ \ of the diferentiable structure on the manifold M
determines a differentiable atlas $A_{Osc^{2}M}$ on $Osc^{2}M$ and therefore
the triple $\left( Osc^{2}M\text{,}\pi ^{2},M\right) $ is a differentiable
bundle.

By $\left( \text{\ref{1.22}}\right) $, a transformation of local coordinates 
$\left( x^{a},y^{\left( 1\right) a},y^{\left( 2\right) a}\right) \rightarrow 
$\newline
$\left( \tilde{x}^{a},\tilde{y}^{\left( 1\right) a},\tilde{y}^{\left(
2\right) a}\right) $ on the manifold $Osc^{2}M$ is given by%
\begin{equation}
\left\{ 
\begin{array}{l}
\tilde{x}^{a}=\tilde{x}^{a}\left( x^{1},...,x^{n}\right) ,~\det \left( 
\dfrac{\partial \tilde{x}^{a}}{\partial x^{b}}\right) \neq 0 \\ 
\\ 
\tilde{y}^{\left( 1\right) a}=\dfrac{\partial \tilde{x}^{a}}{\partial x^{b}}%
y^{\left( 1\right) b} \\ 
\\ 
2\tilde{y}^{\left( 2\right) a}=\dfrac{\partial \tilde{y}^{\left( 1\right) a}%
}{\partial x^{b}}y^{\left( 1\right) b}+2\dfrac{\partial \tilde{y}^{\left(
1\right) a}}{\partial y^{\left( 1\right) b}}y^{\left( 2\right) b}.%
\end{array}%
\right.  \tag{1.3}  \label{1.33}
\end{equation}%
One can see that $Osc^{2}M$ is of dimension 3n.

Let us consider the 2-\textbf{tangent structure} $\mathbb{J}$ on $Osc^{2}M$%
\begin{equation*}
\mathbb{J}\left( \dfrac{\partial }{\partial x^{a}}\right) =\dfrac{\partial }{%
\partial y^{\left( 1\right) a}},\quad \mathbb{J}\left( \dfrac{\partial }{%
\partial y^{\left( 1\right) a}}\right) =\dfrac{\partial }{\partial y^{\left(
2\right) a}},\quad \mathbb{J}\left( \dfrac{\partial }{\partial y^{\left(
2\right) a}}\right) =0
\end{equation*}%
where $\left( \dfrac{\partial }{\partial x^{a}}\mid _{u},\dfrac{\partial }{%
\partial y^{\left( 1\right) a}}\mid _{u},\dfrac{\partial }{\partial
y^{\left( 2\right) a}}\mid _{u}\right) $ is the natural basis of tangent
space \linebreak $T_{u}Osc^{2}M$, $u\in Osc^{2}M.$ If N is a nonlinear
connection on $Osc^{2}M,$ then $N_{0}=N,~\mathbb{J}\left( N_{0}\right) =N_{1}
$ are two distributions geometrically defined on $Osc^{2}M$, all of
dimension n. Let us consider the distributions $V_{2}$ on $Osc^{2}M~$locally
generated by the vector fields $\left\{ \dfrac{\partial }{\partial y^{\left(
2\right) a}}\right\} .$ Consequently, the tangent bundle to $Osc^{2}M$ at
the point $u\in Osc^{2}M$ is given by a direct sum of the vector space:

\begin{equation}
T_{u}Osc^{2}M=N_{0}\left( u\right) \oplus N_{1}\left( u\right) \oplus
V_{2}\left( u\right) ,\forall u\in Osc^{2}M.  \tag{1.4}  \label{4.10}
\end{equation}

We consider $\left\{ \dfrac{\delta }{\delta x^{a}},\dfrac{\delta }{\delta
y^{1\left( a\right) }},\dfrac{\delta }{\delta y^{\left( 2\right) \left(
a\right) }}\right\} $ the adapted basis to the decomposition $\left( \text{%
\ref{4.10}}\right) $ and its dual basis denoted by $\left( dx^{a},\delta
y^{\left( 1\right) a},\delta y^{\left( 2\right) a}\right) ,$ where%
\begin{equation}
\left\{ 
\begin{array}{l}
\dfrac{\delta }{\delta x^{a}}=\dfrac{\partial }{\partial x^{a}}-\underset{%
\left( 1\right) }{N}\overset{}{^{b}}_{a}\dfrac{\delta }{\delta y^{\left(
1\right) b}}-\underset{\left( 2\right) }{N}\overset{}{^{b}}_{a}\dfrac{%
\partial }{\partial y^{\left( 2\right) b}} \\ 
\\ 
\dfrac{\delta }{\delta y^{\left( 1\right) a}}=\qquad \quad \quad \dfrac{%
\partial }{\partial y^{\left( 1\right) a}}-\underset{\left( 1\right) }{N}%
\overset{}{^{b}}_{a}\dfrac{\partial }{\partial y^{\left( 2\right) b}} \\ 
\\ 
\dfrac{\delta }{\delta y^{\left( 2\right) a}}=\qquad \qquad \qquad \qquad
\qquad \dfrac{\partial }{\partial y^{\left( 2\right) a}}%
\end{array}%
\right.  \tag{1.5}  \label{bzad}
\end{equation}

and%
\begin{equation}
\left\{ 
\begin{array}{l}
dx^{a}=\qquad \qquad \qquad \qquad \qquad \qquad dx^{a} \\ 
\\ 
\delta y^{\left( 1\right) a}=\qquad \qquad \quad dy^{\left( 1\right) a}+%
\underset{\left( 1\right) }{M}\overset{}{^{a}}_{b}dx^{b} \\ 
\\ 
\delta y^{\left( 2\right) a}=dy^{\left( 2\right) a}+\underset{\left(
1\right) }{M}\overset{}{^{a}}_{b}\delta y^{b}+\underset{\left( 2\right) }{M}%
\overset{}{^{a}}_{b}\delta y^{\left( 2\right) b}.%
\end{array}%
\right.  \tag{1.6}  \label{1.77}
\end{equation}

\textbf{Definition 1.1} A linear connection $D$ on $Osc^{2}M$\textit{\ }is
called \textbf{N-linear connection} if it preserves by parallelism the
horizontal and vertical distributions $N_{0},N_{1}$ and $V_{2}$ on $%
Osc^{2}M. $

Any N-linear connection $D$ can be represented by an unique system of
functions $D\Gamma \left( N\right) =\left( \underset{\left( i0\right) }{L}%
\overset{}{_{bd}^{a}},\underset{\left( i1\right) }{C}\overset{}{_{bd}^{a}},%
\underset{\left( i2\right) }{C}\overset{}{_{bd}^{a}}\right) ,\left(
i=0,1,2\right) .$ These functions are called \textbf{the coefficients} of
the N-linear connection D.

If on the manifold $Osc^{2}M$ is given a N-linear connection D then there
exists a $h_{i}$-$,v_{1i}$- and $v_{2i}$-\textbf{covariant derivatives} in
local adapted basis $\left( i=0,1,2\right) .$

Any d-tensor $T$, of type $\left( r,s\right) $ can be represented in the
adapted basis and its dual basis in the form%
\begin{equation*}
T=T_{b_{1}...b_{s}}^{a_{1}...a_{r}}\delta _{a_{1}}\otimes ...\otimes \dot{%
\partial}_{2a_{r}}\otimes dx^{b_{1}}\otimes ...\otimes \delta y^{\left(
2\right) b_{s}}.
\end{equation*}%
and we have%
\begin{equation*}
\begin{array}{l}
T_{b_{1}...b_{s}\mid _{id}}^{a_{1}...a_{r}}=\delta
_{a}T_{b_{1}...b_{s}}^{a_{1}...a_{r}}+\underset{\left( i0\right) }{L}\overset%
{}{_{cd}^{a_{1}}}T_{b_{1}...b_{s}}^{ca_{2}...a_{r}}+...+ \\ 
+\underset{\left( i0\right) }{L}\overset{}{_{cd}^{a_{r}}}%
T_{b_{1}...b_{s}}^{a_{1}...a_{r-1}c}-\underset{\left( i0\right) }{L}\overset{%
}{_{b_{1}d}^{c}}T_{cb_{2}...b_{s}}^{a_{1}...a_{r}}-...-\underset{\left(
i0\right) }{L}\overset{}{_{b_{s}d}^{c}}T_{cb_{2}...b_{s-1}c}^{a_{1}...a_{r}},
\\ 
\\ 
T_{b_{1}...b_{s}}^{a_{1}...a_{r}}\overset{\left( 1\right) }{\mid }%
_{id}=\delta _{1a}T_{b_{1}...b_{s}}^{a_{1}...a_{r}}+\underset{\left(
i1\right) }{C}\overset{}{_{cd}^{a_{1}}}%
T_{b_{1}...b_{s}}^{ca_{2}...a_{r}}+...+ \\ 
+\underset{\left( i1\right) }{C}\overset{}{_{cd}^{a_{r}}}%
T_{b_{1}...b_{s}}^{a_{1}...a_{r-1}c}-\underset{\left( i1\right) }{C}\overset{%
}{_{b_{1}d}^{c}}T_{cb_{2}...b_{s}}^{a_{1}...a_{r}}-...-\underset{\left(
i1\right) }{C}\overset{}{_{b_{s}d}^{c}}T_{cb_{2}...b_{s-1}c}^{a_{1}...a_{r}},%
\end{array}%
\end{equation*}%
\begin{equation*}
\begin{array}{l}
T_{b_{1}...b_{s}}^{a_{1}...a_{r}}\overset{\left( 2\right) }{\mid }%
_{id}=\delta _{2a}T_{b_{1}...b_{s}}^{a_{1}...a_{r}}+\underset{\left(
i2\right) }{C}\overset{}{_{cd}^{a_{1}}}%
T_{b_{1}...b_{s}}^{ca_{2}...a_{r}}+...+ \\ 
+\underset{\left( i2\right) }{C}\overset{}{_{cd}^{a_{r}}}%
T_{b_{1}...b_{s}}^{a_{1}...a_{r-1}c}-\underset{\left( i2\right) }{C}\overset{%
}{_{b_{1}d}^{c}}T_{cb_{2}...b_{s}}^{a_{1}...a_{r}}-...-\underset{\left(
i2\right) }{C}\overset{}{_{b_{s}d}^{c}}T_{cb_{2}...b_{s-1}c}^{a_{1}...a_{r}},
\\ 
\\ 
\left( \delta _{1a}=\dfrac{\delta }{\delta y^{\left( 1\right) a}},\delta
_{2a}=\dfrac{\delta }{\delta y^{\left( 2\right) a}};i=0,1,2\right) .%
\end{array}%
\end{equation*}

The operators "$\mid _{id}$" ,"$\overset{\left( 1\right) }{\mid }_{id}$" and
"$\overset{\left( 2\right) }{\mid }_{id}$" are called the h$_{i}$-$,$v$_{1i}$%
- and v$_{2i}$-\textbf{co\-va\-riant derivatives} with respect to $D\Gamma
\left( N\right) .$

\textbf{Definition 1.2 }A\textbf{\ metric structure }on the manifold $%
Osc^{2}\left( M\right) $ is a symmetric covariant tensor field $\mathbb{G}$
of the type $\left( 0,2\right) $ which is non degenerate at each point $u\in
Osc^{2}\left( M\right) $ and of constant signature on $Osc^{2}\left(
M\right) .$

Locally, a metric structure looks as follows:%
\begin{equation*}
\mathbb{G=}\underset{\left( 0\right) }{g}\overset{}{_{ab}}dx^{a}\otimes
dx^{b}+\underset{\left( 1\right) }{g}\overset{}{_{ab}}\delta y^{\left(
1\right) a}\otimes \delta y^{\left( 1\right) b}+\underset{\left( 2\right) }{g%
}\overset{}{_{ab}}\delta y^{\left( 2\right) a}\otimes \delta y^{\left(
2\right) b},
\end{equation*}%
where%
\begin{equation*}
rank\left\Vert \underset{\left( i\right) }{g}\overset{}{_{ab}}\right\Vert
=n,\left( i=0,1,2\right) .
\end{equation*}

\textbf{Definition 1.3 }A N-linear connection $D$ on $Osc^{2}M$ endowed with
a metric structure $\mathbb{G}$ is said to be a \textbf{metric N-linear
connection} if $D_{X}\mathbb{G}=0$ for every $X\in \mathcal{X}\left(
Osc^{2}M\right) .$

\section{Submanifolds in the manifold $Osc^{2}M$}

Let $M$ be a $C^{\infty }$ real, n-dimensional manifold and $\check{M}$ be a
real, m-dimensional manifold, immersed in $M$ through the immersion $i:%
\check{M}\rightarrow M$ . Localy, $i$ can be given in the form%
\begin{equation}
\begin{array}{lcl}
x^{a}=x^{a}\left( u^{1},...,u^{m}\right) , &  & rank\left\Vert \dfrac{%
\partial x^{a}}{\partial u^{\alpha }}\right\Vert =m.%
\end{array}
\tag{2.1}  \label{1.1}
\end{equation}%
The indices $a,b,c,$....run over the set $\left\{ 1,...,n\right\} $ and $%
\alpha ,\beta ,\gamma ,...$ run on the set $\left\{ 1,...,m\right\} .$ We
assume $1<m<n$. If $i$ is an embedding, then we identify $\check{M}$ to $%
i\left( \check{M}\right) $ and say that $\check{M}$ is a \textit{submanifold}
of the manifold M. Therefore $\left( \text{\ref{1.1}}\right) $ will be
called the \textit{parametric equations} \textit{of the submanifold} $\check{%
M}$ in the manifold $M$.

The embedding $i:\check{M}\rightarrow M$ determines an immersion $%
Osc^{2}i:Osc^{2}\check{M}\rightarrow Osc^{2}M$, defined by the covariant
functor $Osc^{2}:Man\rightarrow Man$ \cite{mi}.

The mapping $Osc^{2}i:Osc^{2}\check{M}\rightarrow Osc^{2}M$ has the
parametric equations:%
\begin{equation}
\left\{ 
\begin{array}{l}
x^{a}=x^{a}\left( u^{1},...,u^{m}\right) ,rank\left\Vert \dfrac{\partial
x^{a}}{\partial u^{\alpha }}\right\Vert =m \\ 
\\ 
y^{\left( 1\right) a}=\dfrac{\partial x^{a}}{\partial u^{\alpha }}v^{\left(
1\right) \alpha } \\ 
\\ 
2y^{\left( 2\right) a}=\dfrac{\partial y^{\left( 1\right) a}}{\partial
u^{\alpha }}v^{\left( 1\right) \alpha }+2\dfrac{\partial y^{\left( 1\right)
a}}{\partial v^{\left( 1\right) \alpha }}v^{\left( 2\right) \alpha }%
\end{array}%
\right.  \tag{2.2}  \label{1.3}
\end{equation}%
where%
\begin{equation}
\left\{ 
\begin{array}{c}
\dfrac{\partial x^{a}}{\partial u^{\alpha }}=\dfrac{\partial y^{\left(
1\right) a}}{\partial v^{\left( 1\right) \alpha }}=\dfrac{\partial y^{\left(
2\right) a}}{\partial v^{\left( 2\right) \alpha }} \\ 
\\ 
\dfrac{\partial y^{\left( 1\right) a}}{\partial u^{\alpha }}=\dfrac{\partial
y^{\left( 2\right) a}}{\partial v^{\left( 1\right) \alpha }}.%
\end{array}%
\right.  \tag{2.3}  \label{1.44}
\end{equation}

The Jacobian matrix of $\left( \text{\ref{1.3}}\right) $ is $J\left(
Osc^{2}i\right) $ and it has the rank equal to $3m.$ So, $Osc^{2}i$ is an
immersion. The differential $i_{\ast }$ of the mapping $Osc^{2}i:Osc^{2}%
\check{M}\rightarrow Osc^{2}M$ leads to the relation between the natural
basis of the modules $\mathcal{X}\left( Osc^{2}\check{M}\right) $ and $%
\mathcal{X}\left( Osc^{2}M\right) $ given by%
\begin{equation*}
i_{\ast }\left\Vert \frac{\partial }{\partial u^{\alpha }}\frac{\partial }{%
\partial v^{\left( 1\right) \alpha }}...\frac{\partial }{\partial v^{\left(
2\right) \alpha }}\right\Vert =\left\Vert \frac{\partial }{\partial x^{a}}%
\frac{\partial }{\partial y^{\left( 1\right) a}}...\frac{\partial }{\partial
y^{\left( 2\right) a}}\right\Vert J\left( Osc^{2}i\right) .
\end{equation*}%
$i_{\ast }$ maps the cotangent space $T^{\ast }\left( Osc^{2}M\right) $ in a
point of $Osc^{2}M$, into the cotangent space $T^{\ast }\left( Osc^{2}\check{%
M}\right) $ in a point of $Osc^{2}\check{M}$ by the rule:%
\begin{equation}
\begin{array}{l}
dx^{a}=\dfrac{\partial x^{a}}{\partial u^{\alpha }}du^{\alpha } \\ 
\\ 
dy^{\left( 1\right) a}=\dfrac{\partial y^{\left( 1\right) a}}{\partial
u^{\alpha }}du^{\alpha }+\dfrac{\partial y^{\left( 1\right) a}}{\partial
v^{\left( 1\right) \alpha }}dv^{\left( 1\right) \alpha } \\ 
\\ 
dy^{\left( 2\right) a}=\dfrac{\partial y^{\left( 2\right) a}}{\partial
u^{\alpha }}du^{\alpha }+\dfrac{\partial y^{\left( 2\right) a}}{\partial
v^{\left( 1\right) \alpha }}dv^{\left( 1\right) \alpha }+\dfrac{\partial
y^{\left( 2\right) a}}{\partial v^{\left( 2\right) \alpha }}dv^{\left(
2\right) \alpha }%
\end{array}
\tag{2.4}  \label{2.6}
\end{equation}%
We used the previous theory for study the induced geometrical object fields
from $Osc^{2}M$ to $Osc^{2}\check{M}.$

Let us consider a Finsler space, $F^{n}=(M,F\left( x,y^{\left( 1\right)
}\right) $ having $g_{ab}=\dfrac{1}{2}\dfrac{\partial F^{2}}{\partial
y^{\left( 1\right) a}y^{\left( 1\right) b}}\linebreak$ as fundamental tensor
field. The restriction $\check{F}$ of the function $F$ to the manifold $%
Osc^{2}\check{M}$ is given by%
\begin{equation*}
\check{F}\left( u,v^{\left( 1\right) }\right) =F\left( x\left( u\right)
,y^{\left( 1\right) }\left( u,v^{\left( 1\right) }\right) \right)
\end{equation*}%
and the pair $\check{F}^{n}=\left( M,\check{F}\right) $ is a Finsler space. $%
\check{F}^{n}$ is called the \textbf{induced Finsler subspaces }of the
Finsler space $F^{n}.$

Next, we consider%
\begin{equation}
\overset{}{B_{\alpha }^{a}}=\dfrac{\partial x^{a}}{\partial u^{\alpha }}. 
\tag{2.5}  \label{2.55}
\end{equation}%
and $\mathbb{G=}g_{ab}dx^{a}\otimes dx^{b}+g_{ab}\delta y^{\left( 1\right)
a}\otimes \delta y^{\left( 1\right) a}+g_{ab}\delta y^{\left( 2\right)
a}\otimes \delta y^{\left( 2\right) a}$ \ the Sasaki prolongation of the
fundamental tensor $g$ along $Osc^{2}M$.

There exist a nonlinear connection on the manifold $Osc^{2}M$ determined
only by $g\overset{}{_{ab}}\left( x,y^{\left( 1\right) }\right) $. The dual
coefficients of this nonlinear connection are \cite{mi}: 
\begin{equation}
\begin{array}{l}
\underset{1}{M}\overset{}{^{a}}_{b}=\dfrac{\partial G^{a}}{\partial
y^{\left( 1\right) b}}, \\ 
\\ 
\underset{2}{M}\overset{}{^{a}}_{b}=\dfrac{1}{2}\left( \Gamma \underset{1}{M}%
\overset{}{^{a}}_{b}+\underset{1}{M}\overset{}{^{a}}_{d}\underset{1}{M}%
\overset{}{^{d}}_{b}\right) ,%
\end{array}
\tag{2.6}  \label{cncvar}
\end{equation}%
\newline
where

\begin{equation*}
\begin{array}{l}
G^{a}=\dfrac{1}{2}\gamma \overset{}{_{bc}^{a}}\left( x,y\right) y^{\left(
1\right) b}y^{\left( 1\right) c}, \\ 
\\ 
\Gamma =y^{\left( 1\right) a}\dfrac{\partial }{\partial x^{a}}+2y^{\left(
2\right) a}\dfrac{\partial }{\partial y^{\left( 1\right) a}},%
\end{array}%
\end{equation*}%
and $\gamma \overset{}{_{bc}^{a}}\left( x,y^{\left( 1\right) }\right) $ are
the Christoffel symbols of the fundamental tensor $g$,%
\begin{equation*}
\gamma \overset{}{_{bc}^{a}}\left( x,y^{\left( 1\right) }\right) =\dfrac{1}{2%
}g^{ad}\left( \frac{\partial g_{dc}}{\partial x^{b}}+\frac{\partial g_{bd}}{%
\partial x^{c}}-\frac{\partial g_{bc}}{\partial x^{d}}\right) .
\end{equation*}

Thus, $\left\{ \overset{}{B_{1}^{a}},\overset{}{B_{2}^{a}},...,\overset{}{%
B_{m}^{a}}\right\} $ are m-linear independent d-vector fields on $Osc^{2}%
\check{M}.$ Also, $\left\{ \overset{}{B_{\alpha }^{1}},\overset{}{B_{\alpha
}^{2}},...,\overset{}{B_{\alpha }^{n}}\right\} $ are d-covector fields, with
respect to the next transformations of coordinates:%
\begin{equation}
\left\{ 
\begin{array}{l}
\bar{u}^{\alpha }=\bar{u}^{\alpha }\left( u^{1},...,u^{m}\right)
,rank\left\Vert \dfrac{\partial \bar{u}^{\alpha }}{\partial u^{\beta }}%
\right\Vert =m \\ 
\\ 
\bar{v}^{\left( 1\right) \alpha }=\dfrac{\partial \bar{u}^{\alpha }}{%
\partial u^{\beta }}v^{\left( 1\right) \beta } \\ 
\\ 
2\bar{v}^{\left( 2\right) \alpha }=\dfrac{\partial \bar{v}^{\left( 1\right)
\alpha }}{\partial u^{\beta }}v^{\left( 1\right) \beta }+2\dfrac{\partial 
\bar{v}^{\left( 1\right) \alpha }}{\partial v^{\left( 1\right) \beta }}%
v^{\left( 2\right) \beta }.%
\end{array}%
\right.  \tag{2.7}  \label{2.3}
\end{equation}

Of course, d-vector fields $\left\{ \overset{}{B_{1}^{a}},...,B_{m}^{a}%
\right\} $ are tangent to the submanifold $\check{M}$.

We say that a d-vector field $\xi ^{a}\left( x,y^{\left( 1\right)
},y^{\left( 2\right) }\right) $ is \textbf{normal} to $Osc^{2}\check{M}~$\
if, on $\check{\pi}^{-1}\left( \check{U}\right) \subset Osc^{2}\check{M},$
we have 
\begin{equation*}
g_{ab}\left( x\left( u\right) ,y^{\left( 1\right) }\left( u,v^{\left(
1\right) }\right) \right) B_{\alpha }^{a}\left( u\right) \cdot \xi
^{b}\left( x\left( u\right) ,y^{\left( 1\right) }\left( u,v^{\left( 1\right)
}\right) ,y^{\left( 2\right) }\left( u,v^{\left( 1\right) },v^{\left(
2\right) }\right) \right) =0.
\end{equation*}%
Consequently, on $\check{\pi}^{-1}\left( \check{U}\right) \subset Osc^{2}%
\check{M}$ there exist $n-m$ unit vector fields $B_{\bar{\alpha}}^{a}$, $%
\left( \bar{\alpha}=1,...,n-m\right) $ normal along $Osc^{2}\check{M},$ and
to eatch other:%
\begin{equation}
g_{ab}B_{\alpha }^{a}B_{\bar{\beta}}^{b}=0,\quad g_{ab}B_{\bar{\alpha}%
}^{a}B_{\bar{\beta}}^{b}=\delta _{\bar{\alpha}\bar{\beta}},\left( \bar{\alpha%
},\bar{\beta}=1,...,n-m\right) .  \tag{2.8}  \label{2.22(2.66)}
\end{equation}%
The system of d-vectors $B_{\bar{\alpha}}^{a}$ $\left( \bar{\alpha}%
=1,...,n-m\right) $ is determined up to orthogonal transformations of the
form%
\begin{equation}
B_{\bar{\alpha}^{\prime }}^{a}=A_{\bar{\alpha}^{\prime }}^{\bar{\beta}}B_{%
\bar{\beta}}^{a},\text{ }\left\Vert A_{\bar{\alpha}^{\prime }}^{\bar{\alpha}%
}\right\Vert \in \mathcal{O}\left( n-m\right) ,  \tag{2.9}  \label{2.7}
\end{equation}%
where $\bar{\alpha},\bar{\beta},...$ run over the set $\left\{
1,2,..,n-m\right\} .$

We say that the system of d-vectors $\left\{ B_{\alpha }^{a},B_{\bar{\alpha}%
}^{a}\right\} $ determines a frame in $Osc^{2}M$ along to $Osc^{2}\check{M}.$

Its dual frame will be denoted by $\left\{ B_{a}^{\alpha }\left( u,v^{\left(
1\right) },v^{\left( 2\right) }\right) ,B_{a}^{\bar{\alpha}}\left(
u,v^{\left( 1\right) },v^{\left( 2\right) }\right) \right\} .$ This is also
defined on an open set $\check{\pi}^{-1}\left( \check{U}\right) \subset
Osc^{2}\check{M},$ $\check{U}$ being a domain of a local chart on the
submanifold $\check{M}.$

The conditions of duality are given by:%
\begin{equation}
B_{\beta }^{a}B_{a}^{\alpha }=\delta _{\beta }^{\alpha },\quad B_{\beta
}^{a}B_{a}^{\bar{\alpha}}=0,\quad B_{a}^{\alpha }B_{\bar{\beta}}^{a}=0,\quad
B_{a}^{\bar{\alpha}}B_{\bar{\beta}}^{a}=\delta _{\bar{\beta}}^{\bar{\alpha}}
\tag{2.10}  \label{2.99}
\end{equation}%
\begin{equation}
\begin{array}{l}
B_{\alpha }^{a}B_{b}^{\alpha }+B_{\bar{\alpha}}^{a}B_{b}^{\bar{\alpha}%
}=\delta _{b}^{a}.%
\end{array}
\tag{2.11}  \label{2.99'}
\end{equation}%
Using $\left( \text{\ref{2.22(2.66)}}\right) $, we deduce:%
\begin{equation}
g_{\alpha \beta }B_{a}^{\alpha }=g_{ab}B_{\beta }^{a},\text{~}\delta _{\bar{%
\alpha}\bar{\beta}}B_{b}^{\bar{\beta}}=g_{ab}B_{\bar{\alpha}}^{a}. 
\tag{2.12}  \label{2.10}
\end{equation}%
So, we can look to the set%
\begin{equation*}
\mathcal{R=}\left\{ \left( u,v^{\left( 1\right) },v^{\left( 2\right)
}\right) ;B_{\alpha }^{a}\left( u\right) ,B_{\bar{\alpha}}^{a}\left(
u,v^{\left( 1\right) },v^{\left( 2\right) }\right) \right\}
\end{equation*}%
$\left( u,v^{\left( 1\right) },v^{\left( 2\right) }\right) \in \check{\pi}%
^{-1}\left( \check{U}\right) ~$as a moving frame. Now, we shall represent in 
$\mathcal{R}$ the d-tensor fields from the space $Osc^{2}M$ , restricted to
the open set $\check{\pi}^{-1}\left( \check{U}\right) .$

\section{Induced nonlinear connections}

Now, let us consider the nonlinear connection $N$ (\ref{cncvar}) on the
manifold $Osc^{2}M$ . Then its dual coefficients $\underset{\left( 1\right) }%
{M}\overset{}{_{b}^{a}},\underset{\left( 2\right) }{M}\overset{}{_{b}^{a}%
\text{ depends only by the metric }g\text{.}}$ We will prove that the
restriction of the of the nonlinear connection N to $Osc^{2}\check{M}$
uniquely determines an induced nonlinear connection \v{N} on $Osc^{2}\check{M%
}.$ Of course, \v{N} is well determined by means of its dual coefficients $%
\left( \underset{\left( 1\right) }{\check{M}}\overset{}{_{\beta }^{\alpha }},%
\underset{\left( 2\right) }{\check{M}}\overset{}{_{\beta }^{\alpha }}\right) 
$ or by means of its adapted cobasis $\left( du^{\alpha },\delta v^{\left(
1\right) \alpha },\delta v^{\left( 2\right) \alpha }\right) .$

\textbf{Definition 3.1 }A non-linear connection $\check{N}$ on $Osc^{2}%
\check{M}$\ is called \textbf{induced} by the nonlinear connection N if we
have%
\begin{equation}
\delta v^{\left( 1\right) \alpha }=B_{a}^{\alpha }\delta y^{\left( 1\right)
a},\quad \delta v^{\left( 2\right) \alpha }=B_{a}^{\alpha }\delta y^{\left(
2\right) a}  \tag{3.1}  \label{3.1}
\end{equation}

\textbf{Proposition 3.1} \ \textit{The dual coefficients of the nonlinear
connection \v{N} are} 
\begin{equation}
\begin{array}{l}
\underset{\left( 1\right) }{\check{M}}\overset{}{_{\beta }^{\alpha }}%
=B_{a}^{\alpha }\left( B_{0\beta }^{a}+\underset{\left( 1\right) }{M}\overset%
{}{_{b}^{a}}B_{\beta }^{b}\right) \\ 
\\ 
\underset{\left( 2\right) }{\check{M}}\overset{}{_{\beta }^{\alpha }}%
=B_{a}^{\alpha }\left( \dfrac{1}{2}\dfrac{\partial B_{\delta \gamma }^{a}}{%
\partial u^{\beta }}v^{\left( 1\right) \delta }v^{\left( 1\right) \gamma
}+B_{\delta \beta }^{a}v^{\left( 2\right) \delta }+\underset{\left( 1\right) 
}{M}\overset{}{_{b}^{a}}B_{0\beta }^{b}+\underset{\left( 2\right) }{M}%
\overset{}{_{b}^{a}}B_{\beta }^{b}\right)%
\end{array}
\tag{3.2}  \label{3.2}
\end{equation}%
\textit{where }$\underset{\left( 1\right) }{M}\overset{}{_{b}^{a}},~\underset%
{\left( 2\right) }{M}\overset{}{_{b}^{_{a}}}$\textit{\ are the dual
coefficients of the non-linear connection N.}

\textbf{Theorem 3.1} \ \textit{The cobasis }$\left( dx^{a},\delta y^{\left(
1\right) a},\delta y^{\left( 2\right) a}\right) $\textit{\ restricted to }$%
Osc^{2}\check{M}$ \textit{is}\newline
\textit{uniquely represented in the moving frame }$\mathcal{R}$\textit{\ in
the following form:}%
\begin{equation}
\left\{ 
\begin{array}{l}
dx^{a}=B_{\beta }^{a}du^{\beta } \\ 
\\ 
\delta y^{\left( 1\right) a}=B_{\alpha }^{a}\delta v^{\left( 1\right) \alpha
}+B_{\bar{\alpha}}^{a}\underset{\left( 1\right) }{K}\overset{}{_{\beta }^{%
\bar{\alpha}}}du^{\beta } \\ 
\\ 
\delta y^{\left( 2\right) a}=B_{\alpha }^{a}\delta v^{\left( 2\right) \alpha
}+B_{\bar{\beta}}^{a}\underset{\left( 1\right) }{K}\overset{}{_{\alpha }^{%
\bar{\beta}}}\delta v^{\left( 1\right) \alpha }+B_{\bar{\beta}}^{a}\underset{%
\left( 2\right) }{K}\overset{}{_{\alpha }^{\bar{\beta}}}du^{\alpha }%
\end{array}%
\right.  \tag{3.3}  \label{3.3}
\end{equation}%
\textit{where}%
\begin{equation}
\begin{array}{l}
\underset{\left( 1\right) }{K}\overset{}{_{\beta }^{\bar{\alpha}}}=B_{a}^{%
\bar{\alpha}}\left( B_{0\beta }^{a}+\underset{\left( 1\right) }{M}\overset{}{%
_{b}^{a}}B_{\beta }^{b}\right) \\ 
\\ 
\underset{\left( 2\right) }{K}\overset{}{_{\beta }^{\bar{\alpha}}}=B_{a}^{%
\bar{\alpha}}\left( \dfrac{1}{2}\dfrac{\partial B_{\delta \gamma }^{a}}{%
\partial u^{\beta }}v^{\left( 1\right) \delta }v^{\left( 1\right) \gamma
}+B_{\delta \beta }^{b}v^{\left( 2\right) \delta }+\underset{\left( 1\right) 
}{M}\overset{}{_{b}^{a}}B_{0\beta }^{b}+\underset{\left( 2\right) }{M}%
\overset{}{_{b}^{a}}B_{\beta }^{b}\right. - \\ 
\\ 
-B_{f}^{\bar{\alpha}}B_{d}^{\gamma }\left( B_{\gamma }^{f}+\underset{\left(
1\right) }{M}\overset{}{_{b}^{f}}B_{\gamma }^{b}\right) \left( B_{0\beta
}^{d}+\underset{\left( 1\right) }{M}\overset{}{_{g}^{d}}B_{\beta }^{g}\right)%
\end{array}
\tag{3.4}  \label{3.4}
\end{equation}%
\ \textit{are mixed d-tensor fields}.

\begin{proof}
\textit{The first relation is obviously. From \ref{1.3} and\ \ref{3.2} we
obtain \ref{3.3}}$.$
\end{proof}

\bigskip Generally, a set of functions $T_{j...\beta ...\bar{\beta}%
}^{i...\alpha ...\bar{\alpha}}\left( u,v^{\left( 1\right) },v^{\left(
2\right) }\right) $ which are d-tensors in the index $a,b$,..., and
d-tensors in the index $\alpha ,\beta ,...$and tensors with respect to the
transformations $\left( \text{\ref{2.7}}\right) $ in the index $\bar{\alpha},%
\bar{\beta},...$ is calld a \textit{mixed} d-tensor field on the manifold $%
Osc^{2}\check{M}.$

\section{The relative covariant derivatives}

We shall construct the operators $\underset{\left( i\right) }{\nabla }$ of
relative (or mixed) covariant derivation in the algebra of mixed d-tensor
fields. It is clear that $\underset{\left( i\right) }{\nabla }$ will be
known if its action of functions and on the vector fields of the form%
\begin{equation}
\begin{array}{l}
X^{a}\left( x\left( u\right) ,y^{\left( 1\right) }\left( u,v^{\left(
1\right) }\right) ,y^{\left( 2\right) }\left( u,v^{\left( 1\right)
},v^{\left( 2\right) }\right) \right) \\ 
\\ 
X^{\alpha }\left( u,v^{\left( 1\right) },v^{\left( 2\right) }\right) ,~X^{%
\bar{\alpha}}\left( u,v^{\left( 1\right) },v^{\left( 2\right) }\right)%
\end{array}
\tag{4.1}  \label{3.5}
\end{equation}%
are known.

Let D the canonical metrical N-linear connection on the manifold $Osc^{2}M$ $%
\left( \text{\cite{at},\cite{At+RUS},\cite{At+RUS2}}\right) $%
\begin{equation}
\begin{array}{lll}
\underset{\left( 00\right) }{\overset{c}{L}}\overset{}{_{bc}^{a}} & = & 
\dfrac{1}{2}g\overset{}{^{ad}}\left( \delta _{b}g\underset{}{_{dc}}+\delta
_{c}g\overset{}{_{bd}}-\delta _{d}g\underset{}{_{bc}}\right) , \\ 
&  &  \\ 
\underset{\left( j0\right) }{\overset{c}{L}}\overset{}{_{bc}^{a}} & = & 
\underset{\left( jj\right) }{B}\overset{}{_{cb}^{a}}+\dfrac{1}{2}g\overset{}{%
^{ad}}\left( \delta _{c}g\overset{}{_{bd}}-\underset{\left( jj\right) }{B}%
\overset{}{_{cb}^{f}}g\overset{}{_{fd}}-\underset{\left( jj\right) }{B}%
\overset{}{_{cd}^{f}}g\overset{}{_{bf}}\right) ,\left( j=1,2\right) \\ 
&  &  \\ 
\underset{\left( k1\right) }{\overset{c}{C}}\overset{}{_{bc}^{a}} & = & 
\dfrac{1}{2}g\overset{}{^{ad}}\delta _{1c}g\overset{}{_{bd}},\left(
k=0,2\right) \\ 
&  &  \\ 
\underset{\left( l2\right) }{\overset{c}{C}}\overset{}{_{bc}^{a}} & = & 
\dfrac{1}{2}g\overset{}{^{ad}}\dot{\partial}_{2c}g\overset{}{_{bd}},\left(
l=0,1\right) , \\ 
&  &  \\ 
\underset{\left( ii\right) }{\overset{c}{C}}\overset{}{_{bc}^{a}} & = & 
\dfrac{1}{2}g\overset{}{^{ad}}\left( \delta _{ib}g\underset{}{_{dc}}+\delta
_{ic}g\overset{}{_{bd}}-\delta _{id}g\underset{}{_{bc}}\right) ,\left(
i=1,2\right) .%
\end{array}
\tag{4.2}  \label{N-con.can}
\end{equation}

\textbf{Definition 4.1 }The \textbf{coupling} of the canonical metrical $N$%
-linear connection $D$\textbf{\ }to the induced nonlinear connection $\check{%
N}$ along $Osc^{2}\check{M}$ is locally given by\ the set of its nine
coefficients $\check{D}\Gamma \left( \check{N}\right) =\left( \underset{%
\left( i0\right) }{\check{L}}\overset{}{_{b\delta }^{a}},\underset{\left(
i1\right) }{\check{C}}\overset{}{_{b\delta }^{a}},\underset{\left( i2\right) 
}{\check{C}}\overset{}{_{b\delta }^{a}}\right) ,\left( i=0,1,2\right) ,$
where

\begin{equation}
\begin{array}{l}
\underset{\left( i0\right) }{\check{L}}\overset{}{_{b\delta }^{a}}=\underset{%
\left( i0\right) }{L}\overset{}{_{bd}^{a}}B_{\delta }^{d}+\underset{\left(
i1\right) }{C}\overset{}{_{bd}^{a}}B_{\bar{\delta}}^{d}\underset{\left(
1\right) }{K}\overset{}{_{\delta }^{\bar{\delta}}}+\underset{\left(
i2\right) }{C}\overset{}{_{b\delta }^{a}}B_{\bar{\delta}}^{d}\underset{%
\left( 2\right) }{K}\overset{}{_{\delta }^{\bar{\delta}}} \\ 
\\ 
\underset{\left( i1\right) }{\check{C}}\overset{}{_{b\delta }^{a}}=\underset{%
\left( i1\right) }{C}\overset{}{_{bd}^{a}}B_{\delta }^{d}+\underset{\left(
i2\right) }{C}\overset{}{_{bd}^{a}}B_{\bar{\delta}}^{d}\underset{\left(
1\right) }{K}\overset{}{_{\delta }^{\bar{\delta}}} \\ 
\\ 
\underset{\left( i2\right) }{\check{C}}\overset{}{_{b\delta }^{a}}=\underset{%
\left( i2\right) }{C}\overset{}{_{bd}^{a}B_{\delta }^{d}}.%
\end{array}%
\left( i=0,1,2\right)  \tag{4.3}  \label{3.8}
\end{equation}%
We \ have the operators $\underset{\left( i\right) }{\check{D}}$\ and $%
\underset{\left( i\right) }{D}$\ $\left( i=0,1,2\right) $ with the property%
\begin{equation}
\underset{\left( i\right) }{\check{D}}X^{a}=\underset{\left( i\right) }{D}%
X^{a}~\text{(modulo \ref{3.3})},  \tag{4.4}  \label{3.6}
\end{equation}%
where%
\begin{equation}
\underset{\left( i\right) }{D}X^{a}=dX^{a}+X^{b}\underset{\left( i\right) }{%
\omega }\overset{}{_{b}^{a}},  \tag{4.5}  \label{3.7'}
\end{equation}%
and%
\begin{equation}
\underset{\left( i\right) }{\check{D}}X^{a}=dX^{a}+X^{b}\underset{\left(
i\right) }{\check{\omega}}\overset{}{_{b}^{a}}.  \tag{4.6}  \label{3.7}
\end{equation}%
$\underset{\left( i\right) }{\omega }\overset{}{_{b}^{a}}$ are the 1-forms
of the canonical metrical $N$-connection $D.$ $\underset{\left( i\right) }{%
\check{\omega}}\overset{}{_{b}^{a}}$ are \textbf{coupling 1-forms }of the
coupling\textbf{\ }$\check{D}$ respectivelly.

\bigskip Of course, we can write $\underset{\left( i\right) }{\check{D}}%
X^{a} $ in the form%
\begin{equation*}
\underset{\left( i\right) }{\check{D}}X^{a}=X^{a}\shortmid _{i\delta
}du^{\delta }+X^{a}\overset{\left( 1\right) }{\mid }_{i\delta }\delta
v^{\left( 1\right) \delta }+X^{a}\overset{\left( 2\right) }{\mid }_{i\delta
}\delta v^{\left( 2\right) \delta }.
\end{equation*}

\textbf{Definition 4.2 }We call the\textbf{\ induced tangent }$N$\textbf{%
-linear connection }on $Osc^{2}\check{M}$\ by the canonical metrical $N$%
-linear connection $D$ the set of its nine coefficients $D^{\top }\Gamma
\left( \check{N}\right) =\left( \underset{\left( i0\right) }{L}\overset{}{%
_{\beta \delta }^{\alpha }},\underset{\left( i1\right) }{C}\overset{}{%
_{\beta \delta }^{\alpha }},\underset{\left( i2\right) }{C}\overset{}{%
_{\beta \delta }^{\alpha }}\right) $ $\left( i=0,1,2\right) ,$ where%
\begin{equation}
\begin{array}{l}
\underset{\left( i0\right) }{L}\overset{}{_{\beta \delta }^{\alpha }}%
=B_{d}^{\alpha }\left( B_{\beta \delta }^{d}+B_{\beta }^{f}\underset{\left(
i0\right) }{\check{L}}\overset{}{_{f\delta }^{d}}\right) \\ 
\\ 
\underset{\left( i1\right) }{C}\overset{}{_{\beta \delta }^{\alpha }}%
=B_{d}^{\alpha }B_{\beta }^{f}\underset{\left( i1\right) }{\check{C}}\overset%
{}{_{f\delta }^{d}} \\ 
\\ 
\underset{\left( i2\right) }{C}\overset{}{_{\beta \delta }^{\alpha }}%
=B_{d}^{\alpha }B_{\beta }^{f}\underset{\left( i2\right) }{\check{C}}\overset%
{}{_{f\delta }^{d}}.%
\end{array}%
\left( i=0,1,2\right)  \tag{4.7}  \label{3.12}
\end{equation}%
We \ have the operators $\underset{\left( i\right) }{D}^{\top }$ with the
properties

\begin{equation}
\begin{array}{c}
\underset{\left( i\right) }{D}^{\top }X^{\alpha }=B_{b}^{\alpha }\underset{%
\left( i\right) }{\check{D}}X^{b},\quad \text{\textit{for} }X\overset{}{^{a}}%
=B_{\gamma }^{a}X^{\gamma }%
\end{array}
\tag{4.8}  \label{3.9}
\end{equation}%
\begin{equation}
\underset{\left( i\right) }{D}^{\top }X^{\alpha }=dX^{\alpha }+X^{\beta }%
\underset{\left( i\right) }{\omega }\overset{}{_{\beta }^{\alpha }}, 
\tag{4.9}  \label{3.10}
\end{equation}%
where $\underset{\left( i\right) }{\omega }\overset{}{_{\beta }^{\alpha }}$
are the\textbf{\ induced tangent connection 1-forms} of $\underset{\left(
i\right) }{D}^{\top }~\left( i=0,1,2\right) .$

As in the case of $\check{D}$ we may write%
\begin{equation*}
\underset{\left( i\right) }{D}^{\top }X^{\alpha }=X_{\mid i\delta }^{\alpha
}du^{\delta }+X^{\alpha }\overset{\left( 1\right) }{\mid }_{i\delta }\delta
v^{\left( 1\right) \delta }+X^{\alpha }\overset{\left( 2\right) }{\mid }%
_{i\delta }\delta v^{\left( 2\right) \delta }.
\end{equation*}

\textbf{Definition 4.3} \textit{\ }We call the \textbf{induced normal }$N$%
\textbf{-linear connection} on $Osc^{2}\check{M}$\ by the canonical metrical 
$N$-linear connection $D$ the set of its nine coefficients $D^{\bot }\Gamma
\left( \check{N}\right) =\left( \underset{\left( i0\right) }{L}\overset{}{_{%
\bar{\beta}\delta }^{\bar{\alpha}}},\underset{\left( i1\right) }{C}\overset{}%
{_{\bar{\beta}\delta }^{\bar{\alpha}}},\underset{\left( i2\right) }{C}%
\overset{}{_{\bar{\beta}\delta }^{\bar{\alpha}}}\right) ,$ where%
\begin{equation}
\begin{array}{l}
\underset{\left( i0\right) }{L}\overset{}{_{\bar{\beta}\delta }^{\bar{\alpha}%
}}=B_{d}^{\bar{\alpha}}\left( \dfrac{\delta B_{\bar{\beta}}^{d}}{\delta
u^{\delta }}+B_{\bar{\beta}}^{f}\underset{\left( i0\right) }{\check{L}}%
\overset{}{_{f\delta }^{d}}\right) \\ 
\\ 
\underset{\left( i1\right) }{C}\overset{}{_{\bar{\beta}\delta }^{\bar{\alpha}%
}}=B_{d}^{\bar{\alpha}}\left( \dfrac{\delta B_{\beta }^{d}}{\delta v^{\left(
1\right) \delta }}+B_{\bar{\beta}}^{f}\underset{\left( i1\right) }{\check{C}}%
\overset{}{_{f\delta }^{d}}\right) \\ 
\\ 
\underset{\left( i2\right) }{C}\overset{}{_{\bar{\beta}\delta }^{\bar{\alpha}%
}}=B_{d}^{\bar{\alpha}}\left( \dfrac{\partial B_{\beta }^{d}}{\partial
v^{\left( 2\right) \delta }}+B_{\bar{\beta}}^{f}\underset{\left( i2\right) }{%
\check{C}}\overset{}{_{f\delta }^{d}}\right) .%
\end{array}%
\left( i=0,1,2\right)  \tag{4.10}  \label{3.16}
\end{equation}%
As before$,$ we have the operators $\underset{\left( i\right) }{D}^{\bot }$
with the properties 
\begin{equation}
\begin{array}{c}
\underset{\left( i\right) }{D}^{\bot }X^{\bar{\alpha}}=B_{b}^{\bar{\alpha}}%
\underset{\left( i\right) }{\check{D}}X^{b},\quad \text{\textit{for} }%
X^{a}=B_{\bar{\gamma}}^{a}X^{\bar{\gamma}}%
\end{array}
\tag{4.11}  \label{3.13}
\end{equation}%
\begin{equation}
\underset{\left( i\right) }{D}^{\bot }X^{\bar{\alpha}}=dX^{\bar{\alpha}}+X^{%
\bar{\beta}}\underset{\left( i\right) }{\omega }\overset{}{_{\bar{\beta}}^{%
\bar{\alpha}}}  \tag{4.12}  \label{3.14}
\end{equation}%
where $\underset{\left( i\right) }{\omega }\overset{}{_{\bar{\beta}}^{\bar{%
\alpha}}}$ are the \textbf{induced normal connection 1-forms} of $\underset{%
\left( i\right) }{D}^{\bot }~\left( i=0,1,2\right) .$

We may set%
\begin{equation*}
\underset{\left( i\right) }{D}^{\perp }X^{\bar{\alpha}}=X_{\mid i\delta }^{%
\bar{\alpha}}du^{\delta }+X^{\bar{\alpha}}\overset{\left( 1\right) }{\mid }%
_{i\delta }\delta v^{\left( 1\right) \delta }+X^{\bar{\alpha}}\overset{%
\left( 2\right) }{\mid }_{i\delta }\delta v^{\left( 2\right) \delta }.
\end{equation*}%
Now, we can define the relative (or mixed) covariant derivatives $\underset{%
\left( i\right) }{\nabla }$ enounced at the begining of this section.

\textbf{Theorem 4.4} \ \textit{A relative (mixed) covariant derivation in
the algebra of mixed d-tensor fields is an operator }$\underset{\left(
i\right) }{\nabla }$\textit{\ for which the following properties hold}:%
\begin{equation*}
\begin{array}{c}
\underset{\left( i\right) }{\nabla }f=df,\quad \forall f\in \mathcal{F}%
\left( Osc^{2}\check{M}\right) ~ \\ 
\\ 
\underset{\left( i\right) }{\nabla }X^{a}=\underset{\left( i\right) }{\check{%
D}}X^{a},\quad \underset{\left( i\right) }{\nabla }X^{\alpha }=\underset{%
\left( i\right) }{D}^{\intercal }X^{\alpha },\quad \underset{\left( i\right) 
}{\nabla }X^{\bar{\alpha}}=\underset{\left( i\right) }{D}^{\bot }X^{\bar{%
\alpha}}~\left( i=0,1,2\right)%
\end{array}%
\end{equation*}%
\textit{The connection 1-forms }$\underset{\left( i\right) }{\check{\omega}}%
\overset{}{_{b}^{a}},\underset{\left( i\right) }{\omega }\overset{}{_{\beta
}^{\alpha }},\underset{\left( i\right) }{\omega }\overset{}{_{\bar{\beta}}^{%
\bar{\alpha}}}$\textit{\ will be called the\textbf{\ connection 1-forms} of }%
$\ \underset{\left( i\right) }{\nabla }.$

\textbf{The Liouville vector fields} of the submanifold $Osc^{2}\check{M}$,
introduce by R. Miron in \cite{mi}, are%
\begin{equation*}
\begin{array}{l}
\underset{1}{\mathbb{C}}=v^{\left( 1\right) \alpha }\dfrac{\partial }{%
\partial v^{\left( 2\right) \alpha }} \\ 
\\ 
\underset{2}{\mathbb{C}}=v^{\left( 1\right) \alpha }\dfrac{\partial }{%
\partial v^{\left( 1\right) \alpha }}+2v^{\left( 2\right) \alpha }\dfrac{%
\partial }{\partial v^{\left( 2\right) \alpha }}.%
\end{array}%
\end{equation*}

If we represent this vector fields in the adapted basis, we get%
\begin{equation*}
\underset{1}{\mathbb{C}}=z^{\left( 1\right) \alpha }\dot{\partial}_{2\alpha
},\underset{2}{\mathbb{C}}=z^{\left( 1\right) \alpha }\delta _{1\alpha
}+2z^{\left( 2\right) \alpha }\dot{\partial}_{2\alpha }
\end{equation*}%
where%
\begin{equation*}
z^{\left( 1\right) \alpha }=v^{\left( 1\right) \alpha },z^{\left( 2\right)
\alpha }=v^{\left( 2\right) \alpha }+\frac{1}{2}\underset{\left( 1\right) }{M%
}\overset{}{^{\alpha }}_{\beta }v^{\left( 1\right) \beta }.
\end{equation*}%
The d -vector fields $z^{\left( 1\right) \alpha }$ and $z^{\left( 2\right)
\alpha }$ are called the \textbf{Liouville d-vector fields} of the
submanifold $Osc^{2}\check{M}.$

The $\left( z^{\left( 1\right) }\right) $- and $\left( z^{\left( 2\right)
}\right) $-\textbf{deflection} tensor fields are:%
\begin{equation}
\begin{array}{lllll}
z^{\left( 1\right) \alpha }\overset{}{_{\mid i\beta }}=\underset{i}{\overset{%
\left( 1\right) }{D}}\overset{}{_{\beta }^{\alpha }}, &  & z^{\left(
1\right) \alpha }\overset{\left( 1\right) }{\mid }_{i\beta }=\underset{i}{%
\overset{\left( 11\right) }{d}}\overset{}{_{\beta }^{\alpha }}, &  & 
z^{\left( 1\right) \alpha }\overset{\left( 2\right) }{\mid }_{i\beta }=%
\underset{i}{\overset{\left( 12\right) }{d}}\overset{}{_{\beta }^{\alpha }},
\\ 
&  &  &  &  \\ 
z^{\left( 2\right) \alpha }\overset{}{_{\mid i\beta }}=\underset{i}{\overset{%
\left( 2\right) }{D}}\overset{}{_{\beta }^{\alpha }}, &  & z^{\left(
2\right) \alpha }\overset{\left( 1\right) }{\mid }_{i\beta }=\underset{i}{%
\overset{\left( 21\right) }{d}}\overset{}{_{\beta }^{\alpha }}, &  & 
z^{\left( 2\right) \alpha }\overset{\left( 2\right) }{\mid }_{i\beta }=%
\underset{i}{\overset{\left( 22\right) }{d}}\overset{}{_{\beta }^{\alpha }}.%
\end{array}
\tag{4.13}  \label{5.1}
\end{equation}

\textbf{Proposition 4.1} The $\left( z^{\left( 1\right) }\right) $\textit{%
-deflection tensor fields have the expressions:}%
\begin{equation}
\begin{array}{l}
\underset{i}{\overset{\left( 1\right) }{D}}\overset{}{_{\beta }^{\alpha }}=-%
\underset{1}{N}\overset{}{^{\alpha }}\overset{}{_{\beta }}+z^{\left(
1\right) \gamma }\underset{\left( i0\right) }{L}\overset{}{_{\gamma \beta
}^{\alpha }}, \\ 
\\ 
\underset{i}{\overset{\left( 11\right) }{d}}\overset{}{_{\beta }^{\alpha }}%
=\delta _{\beta }^{\alpha }+z^{\left( 1\right) \gamma }\underset{\left(
i1\right) }{C}\overset{}{_{\gamma \beta }^{\alpha }}, \\ 
\\ 
\underset{i}{\overset{\left( 12\right) }{d}}\overset{}{_{\beta }^{\alpha }}%
=z^{\left( 1\right) \gamma }\underset{\left( i2\right) }{C}\overset{}{%
_{\gamma \beta }^{\alpha }},\text{ }\left( i=0,1,2\right) .%
\end{array}
\tag{4.14}  \label{z^1defl.tens}
\end{equation}%
\ 

Indeed, if we take%
\begin{equation*}
\begin{array}{l}
z^{\left( 1\right) \alpha }\overset{}{_{\mid i\beta }}=\delta _{\beta
}z^{\left( 1\right) \alpha }+z^{\left( 1\right) \gamma }\underset{\left(
i0\right) }{L}\overset{}{_{\gamma \beta }^{\alpha }}, \\ 
\\ 
z^{\left( 1\right) \alpha }\overset{\left( j\right) }{\mid }_{i\beta
}=\delta _{j\beta }z^{\left( 1\right) \alpha }+z^{\left( 1\right) \gamma }%
\underset{\left( ij\right) }{C}\overset{}{_{\gamma \beta }^{\alpha }},\text{ 
}\left( i=0,1,2;j=1,2;\delta _{2\beta }=\dot{\partial}_{2\beta }\right)%
\end{array}%
\end{equation*}%
we find this formulae.

\textbf{Proposition 4.2} The $\left( z^{\left( 2\right) }\right) $\textit{%
-deflection tensor fields are given by}%
\begin{equation}
\begin{array}{l}
\underset{i}{\overset{\left( 2\right) }{D}}\overset{}{_{\beta }^{\alpha }}=%
\dfrac{1}{2}\left( \underset{2}{N}\overset{}{^{\alpha }}\overset{}{_{\beta }}%
+\underset{2}{M}\overset{}{^{\alpha }}\overset{}{_{\beta }}\right) +\dfrac{1%
}{2}z^{\left( 1\right) \gamma }\delta _{\beta }\underset{1}{N}\overset{}{%
^{\alpha }}\overset{}{_{\gamma }}+z^{\left( 2\right) \gamma }\underset{%
\left( i0\right) }{L}\overset{}{_{\gamma \beta }^{\alpha }}, \\ 
\\ 
\underset{i}{\overset{\left( 21\right) }{d}}\overset{}{_{\beta }^{\alpha }}=-%
\dfrac{1}{2}\left( 2\underset{2}{N}\overset{}{^{\alpha }}\overset{}{_{\beta }%
}-\underset{1}{N}\overset{}{^{\alpha }}\overset{}{_{\beta }}\right) +\dfrac{1%
}{2}z^{\left( 1\right) \gamma }\underset{11}{B}\overset{}{_{\gamma \beta
}^{\alpha }}+z^{\left( 2\right) \gamma }\underset{\left( i1\right) }{C}%
\overset{}{_{\gamma \beta }^{\alpha }}, \\ 
\\ 
\underset{i}{\overset{\left( 22\right) }{d}}\overset{}{_{\beta }^{\alpha }}%
=\delta _{\beta }^{\alpha }+\dfrac{1}{2}z^{\left( 1\right) \gamma }\underset{%
12}{B}\overset{}{_{\gamma \beta }^{\alpha }}+z^{\left( 2\right) \gamma }%
\underset{\left( i2\right) }{C}\overset{}{_{\gamma \beta }^{\alpha }},\text{ 
}\left( i=0,1,2\right) .%
\end{array}
\tag{4.15}  \label{z^2defl.tens}
\end{equation}

\section{\protect\bigskip The Ricci identities}

Let $\check{D}\Gamma \left( \check{N}\right) =\left( \underset{\left(
i0\right) }{\check{L}}\overset{}{_{b\delta }^{a}},\underset{\left( i1\right) 
}{\check{C}}\overset{}{_{b\delta }^{a}},\underset{\left( i2\right) }{\check{C%
}}\overset{}{_{b\delta }^{a}}\right) $ the coupling\textbf{\ \ }of the
canonical metrical $N$-linear connection $D$ (\ref{N-con.can}) to the
induced nonlinear connection $\check{N}$ along to the manifold $Osc^{2}%
\check{M}$, $D^{\top }\Gamma \left( \check{N}\right) =\left( \underset{%
\left( i0\right) }{L}\overset{}{_{\beta \delta }^{\alpha }},\underset{\left(
i1\right) }{C}\overset{}{_{\beta \delta }^{\alpha }},\underset{\left(
i2\right) }{C}\overset{}{_{\beta \delta }^{\alpha }}\right) $ and $D^{\bot
}\Gamma \left( \check{N}\right) =\linebreak \left( \underset{\left(
i0\right) }{L}\overset{}{_{\bar{\beta}\delta }^{\bar{\alpha}}},\underset{%
\left( i1\right) }{C}\overset{}{_{\bar{\beta}\delta }^{\bar{\alpha}}},%
\underset{\left( i2\right) }{C}\overset{}{_{\bar{\beta}\delta }^{\bar{\alpha}%
}}\right) $ $\left( i=0,1,2\right) $ the induced tangent $N$-linear
connection and the induced normal $N$-linear connection on $Osc^{2}\check{M}$%
,\ respectivelly.

\textbf{Teorem\u{a} 5.1} For any d-vector fields $X^{\alpha },$ the
following Ricci identities hold:%
\begin{equation}
\begin{array}{lll}
X^{\alpha }\overset{}{_{\mid _{i\beta }\mid _{i\gamma }}}-X^{\alpha }\overset%
{}{_{\mid _{i\beta }\mid _{i\gamma }}} & = & X^{\delta }\underset{\left(
0i\right) }{R}\underset{}{_{\delta }}\overset{}{^{\alpha }}\underset{}{%
_{\beta }}\underset{}{_{\gamma }}-\underset{\left( 0\right) }{\overset{%
\left( i\right) }{T}}\overset{}{_{\beta \gamma }^{\sigma }}X^{\alpha }%
\overset{}{_{\mid _{i\sigma }}}-\underset{\left( 01\right) }{R}\overset{}{%
_{\beta \gamma }^{\sigma }}X^{\alpha }\overset{\left( 1\right) }{\mid }%
_{i\sigma }- \\ 
&  &  \\ 
&  & -\underset{\left( 02\right) }{R}\overset{}{_{\beta \gamma }^{\sigma }}%
X^{\alpha }\overset{\left( 2\right) }{\mid }_{i\sigma }, \\ 
&  &  \\ 
X^{\alpha }\overset{}{_{\mid _{i\beta }}}\overset{\left( 1\right) }{\mid }%
_{i\gamma }-X^{\alpha }\overset{\left( 1\right) }{\mid }_{i\gamma }\overset{}%
{_{\mid _{i\beta }}} & = & X^{\delta }\underset{\left( 1i\right) }{P}%
\underset{}{_{\delta }}\overset{}{^{\alpha }}\underset{}{_{\beta }}\underset{%
}{_{\gamma }}-\underset{\left( i1\right) }{C}\overset{}{_{\beta \gamma
}^{\sigma }}X^{\alpha }\overset{}{_{\mid _{i\sigma }}}-\underset{\left(
11\right) }{\overset{\left( i\right) }{P}}\overset{}{_{\beta \gamma
}^{\sigma }}X^{\alpha }\overset{\left( 1\right) }{\mid }_{i\sigma } \\ 
&  &  \\ 
&  & -\underset{\left( 12\right) }{P}\overset{}{_{\beta \gamma }^{\sigma }}%
X^{\alpha }\overset{\left( 2\right) }{\mid }_{i\sigma }, \\ 
&  &  \\ 
X^{\alpha }\overset{}{_{\mid _{i\beta }}}\overset{\left( 2\right) }{\mid }%
_{i\gamma }-X^{\alpha }\overset{\left( 2\right) }{\mid }_{i\gamma }\overset{}%
{_{\mid _{i\beta }}} & = & X^{\delta }\underset{\left( 2i\right) }{P}%
\underset{}{_{\delta }}\overset{}{^{\alpha }}\underset{}{_{\beta }}\underset{%
}{_{\gamma }}-\underset{\left( i2\right) }{C}\overset{}{_{\beta \gamma
}^{\sigma }}X^{\alpha }\overset{}{_{\mid _{i\sigma }}}-\underset{\left(
21\right) }{P}\overset{}{_{\beta \gamma }^{\sigma }}X^{\alpha }\overset{%
\left( 1\right) }{\mid }_{i\sigma } \\ 
&  &  \\ 
&  & -\underset{\left( 22\right) }{\overset{\left( i\right) }{P}}\overset{}{%
_{\beta \gamma }^{\sigma }}X^{\alpha }\overset{\left( 2\right) }{\mid }%
_{i\sigma },%
\end{array}
\tag{5.1}  \label{5.3}
\end{equation}

\begin{equation*}
\begin{array}{lll}
X^{\alpha }\overset{\left( 1\right) }{\mid }_{i\beta }\overset{\left(
2\right) }{\mid }_{i\gamma }-X^{\alpha }\overset{\left( 2\right) }{\mid }%
_{i\gamma }\overset{\left( 1\right) }{\mid }_{i\beta } & = & X^{\delta }%
\underset{\left( 21\right) }{Q}\underset{}{_{\delta }}\overset{}{^{\alpha }}%
\underset{}{_{\beta }}\underset{}{_{\gamma }}-\underset{\left( i2\right) }{C}%
\overset{}{_{\beta \gamma }^{\sigma }}X^{\alpha }\overset{\left( 1\right) }{%
\mid }_{i\sigma }- \\ 
&  &  \\ 
&  & -\underset{\left( 22\right) }{\overset{\left( i\right) }{Q}}\overset{}{%
_{\beta \gamma }^{\sigma }}X^{\alpha }\overset{\left( 2\right) }{\mid }%
_{i\sigma }, \\ 
&  &  \\ 
X^{\alpha }\overset{\left( j\right) }{\mid }_{i\beta }\overset{\left(
j\right) }{\mid }_{i\gamma }-X^{\alpha }\overset{\left( j\right) }{\mid }%
_{i\gamma }\overset{\left( j\right) }{\mid }_{i\beta } & = & X^{\delta }%
\underset{\left( ji\right) }{S}\underset{}{_{\delta }}\overset{}{^{\alpha }}%
\underset{}{_{\beta }}\underset{}{_{\gamma }}-\underset{\left( j\right) }{%
\overset{\left( i\right) }{S}}\overset{}{_{\beta \gamma }^{\sigma }}%
X^{\alpha }\overset{\left( 1\right) }{\mid }_{i\sigma }- \\ 
&  &  \\ 
&  & -\underset{\left( j2\right) }{R}\overset{}{_{\beta \gamma }^{\sigma }}%
X^{\alpha }\overset{\left( 2\right) }{\mid }_{i\sigma }%
\end{array}%
\end{equation*}%
where $\underset{\left( 22\right) }{R}\overset{}{_{\beta \gamma }^{\alpha }}%
=0,\left( i=0,1,2,j=1,2\right) .$

The Ricci identities (\ref{5.3}) applied to the Liouville d-vector fields $%
z^{\left( 1\right) \alpha }$ and $z^{\left( 2\right) \alpha }$ lead to the
some fundamental identities.

\textbf{Theorem 5.2 }\textit{The deflection tensor fields satisfy the
following identities:}%
\begin{equation}
\begin{array}{lll}
\underset{i}{\overset{\left( j\right) }{D}}\overset{}{_{\beta }^{\alpha }}%
\overset{}{_{\mid _{i\gamma }}}-\underset{i}{\overset{\left( j\right) }{D}}%
\overset{}{_{\gamma }^{\alpha }}\overset{}{_{\mid _{i\beta }}} & = & 
z^{\left( j\right) \delta }\underset{\left( 0i\right) }{R}\underset{}{%
_{\delta }}\overset{}{^{\alpha }}\underset{}{_{\beta }}\underset{}{_{\gamma }%
}-\underset{i}{\overset{\left( j\right) }{D}}\overset{}{_{\delta }^{\alpha }}%
\underset{\left( 0\right) }{\overset{\left( i\right) }{T}}\overset{}{_{\beta
\gamma }^{\delta }}- \\ 
&  &  \\ 
&  & -\underset{i}{\overset{\left( j1\right) }{d}}\overset{}{_{\delta
}^{\alpha }}\underset{\left( 01\right) }{R}\overset{}{_{\beta \gamma
}^{\delta }}-\underset{i}{\overset{\left( j2\right) }{d}}\overset{}{_{\delta
}^{\alpha }}\underset{\left( 02\right) }{R}\overset{}{_{\beta \gamma
}^{\delta }}, \\ 
&  &  \\ 
\underset{i}{\overset{\left( j\right) }{D}}\overset{}{_{\beta }^{\alpha }}%
\overset{\left( 1\right) }{\mid }_{i\gamma }-\underset{i}{\overset{\left(
j1\right) }{d}}\overset{}{_{\gamma }^{\alpha }}\overset{}{_{\mid _{i\beta }}}
& = & z^{\left( j\right) \delta }\underset{\left( 1i\right) }{P}\underset{}{%
_{\delta }}\overset{}{^{\alpha }}\underset{}{_{\beta }}\underset{}{_{\gamma }%
}-\underset{i}{\overset{\left( j\right) }{D}}\overset{}{_{\delta }^{\alpha }}%
\underset{\left( i1\right) }{C}\overset{}{_{\beta \gamma }^{\delta }}- \\ 
&  &  \\ 
&  & -\underset{i}{\overset{\left( j1\right) }{d}}\overset{}{_{\delta
}^{\alpha }}\underset{\left( 11\right) }{\overset{\left( i\right) }{P}}%
\overset{}{_{\beta \gamma }^{\sigma }}-\underset{i}{\overset{\left(
j2\right) }{d}}\overset{}{_{\delta }^{\alpha }}\underset{\left( 12\right) }{P%
}\overset{}{_{\beta \gamma }^{\delta }}, \\ 
&  &  \\ 
\underset{i}{\overset{\left( j\right) }{D}}\overset{}{_{\beta }^{\alpha }}%
\overset{\left( 2\right) }{\mid }_{i\gamma }-\underset{i}{\overset{\left(
j2\right) }{d}}\overset{}{_{\gamma }^{\alpha }}\overset{}{_{\mid _{i\beta }}}
& = & z^{\left( j\right) \delta }\underset{\left( 2i\right) }{P}\underset{}{%
_{\delta }}\overset{}{^{\alpha }}\underset{}{_{\beta }}\underset{}{_{\gamma }%
}-\underset{i}{\overset{\left( j\right) }{D}}\overset{}{_{\delta }^{\alpha }}%
\underset{\left( i2\right) }{C}\overset{}{_{\beta \gamma }^{\sigma }}- \\ 
&  &  \\ 
&  & -\underset{i}{\overset{\left( j1\right) }{d}}\overset{}{_{\delta
}^{\alpha }}\underset{\left( 21\right) }{P}\overset{}{_{\beta \gamma
}^{\delta }}-\underset{i}{\overset{\left( j2\right) }{d}}\overset{}{_{\delta
}^{\alpha }}\underset{\left( 22\right) }{\overset{\left( i\right) }{P}}%
\overset{}{_{\beta \gamma }^{\sigma }}, \\ 
&  &  \\ 
\underset{i}{\overset{\left( j1\right) }{d}}\overset{}{_{\beta }^{\alpha }}%
\overset{\left( 2\right) }{\mid }_{i\gamma }-\underset{i}{\overset{\left(
j2\right) }{d}}\overset{}{_{\gamma }^{\alpha }}\overset{\left( 1\right) }{%
\mid }_{i\beta } & = & z^{\left( j\right) \delta }\underset{\left( 2i\right) 
}{Q}\underset{}{_{\delta }}\overset{}{^{\alpha }}\underset{}{_{\beta }}%
\underset{}{_{\gamma }}-X^{\alpha }\overset{\left( 1\right) }{\mid }%
_{i\sigma }- \\ 
&  &  \\ 
&  & -\underset{i}{\overset{\left( j1\right) }{d}}\overset{}{_{\delta
}^{\alpha }}\underset{\left( i2\right) }{C}\overset{}{_{\beta \gamma
}^{\delta }}-\underset{i}{\overset{\left( j2\right) }{d}}\overset{}{_{\delta
}^{\alpha }}\underset{\left( 22\right) }{\overset{\left( i\right) }{Q}}%
\overset{}{_{\beta \gamma }^{\delta }},%
\end{array}
\tag{5.2}  \label{5.4}
\end{equation}%
\begin{equation*}
\begin{array}{lll}
\underset{i}{\overset{\left( ji\right) }{d}}\overset{}{_{\beta }^{\alpha }}%
\overset{\left( l\right) }{\mid }_{i\gamma }-\underset{i}{\overset{\left(
ji\right) }{d}}\overset{}{_{\gamma }^{\alpha }}\overset{\left( l\right) }{%
\mid }_{i\beta } & = & z^{\left( j\right) \delta }\underset{\left( li\right) 
}{S}\underset{}{_{\delta }}\overset{}{^{\alpha }}\underset{}{_{\beta }}%
\underset{}{_{\gamma }}- \\ 
&  &  \\ 
&  & -\underset{i}{\overset{\left( jl\right) }{d}}\overset{}{_{\delta
}^{\alpha }}\underset{\left( i\right) }{\overset{\left( l\right) }{S}}%
\overset{}{_{\beta \gamma }^{\delta }}-\underset{i}{\overset{\left(
j2\right) }{d}}\overset{}{_{\delta }^{\alpha }}\underset{\left( l2\right) }{R%
}\overset{}{_{\beta \gamma }^{\delta }},%
\end{array}%
\end{equation*}%
$\left( i=0,1,2;j,l=1,2;\underset{\left( 22\right) }{R}\overset{}{_{\beta
\gamma }^{\alpha }}=0.\right) $

Also, if the $\left( z^{\left( 1\right) }\right) $-and $\left( z^{\left(
2\right) }\right) $-deflection tensors have the following particular form%
\begin{equation}
\begin{array}{lllll}
\underset{i}{\overset{\left( 1\right) }{D}}\overset{}{_{\beta }^{\alpha }}=0,
&  & \underset{i}{\overset{\left( 11\right) }{d}}\overset{}{_{\beta
}^{\alpha }}=\delta _{\beta }^{\alpha }, &  & \underset{i}{\overset{\left(
12\right) }{d}}\overset{}{_{\beta }^{\alpha }}=0 \\ 
&  &  &  &  \\ 
\underset{i}{\overset{\left( 2\right) }{D}}\overset{}{_{\beta }^{\alpha }}=0,
&  & \underset{i}{\overset{\left( 21\right) }{d}}\overset{}{_{\beta
}^{\alpha }}=0, &  & \underset{i}{\overset{\left( 22\right) }{d}}\overset{}{%
_{\beta }^{\alpha }}=\delta _{\beta }^{\alpha }%
\end{array}
\tag{5.3}  \label{5.5}
\end{equation}%
then, the fundamental identities from (\ref{5.4}) are very important,
especially for applications.

\textbf{Proposition 5.2 }If the deflection tensors are given by (\ref{5.5}),
then the following identities hold: 
\begin{equation}
\begin{array}{lllll}
z^{\left( j\right) \delta }\underset{\left( 0i\right) }{R}\underset{}{%
_{\delta }}\overset{}{^{\alpha }}\underset{}{_{\beta }}\underset{}{_{\gamma }%
}=\underset{\left( 0j\right) }{R}\overset{}{_{\beta \gamma }^{\alpha }}, & 
& z^{\left( 1\right) \delta }\underset{\left( 2i\right) }{P}\underset{}{%
_{\delta }}\overset{}{^{\alpha }}\underset{}{_{\beta }}\underset{}{_{\gamma }%
}=\underset{\left( 21\right) }{P}\overset{}{_{\beta \gamma }^{\alpha }}, & 
& z^{\left( 2\right) \delta }\underset{\left( 1i\right) }{P}\underset{}{%
_{\delta }}\overset{}{^{\alpha }}\underset{}{_{\beta }}\underset{}{_{\gamma }%
}=\underset{\left( 12\right) }{P}\overset{}{_{\beta \gamma }^{\alpha }}, \\ 
&  &  &  &  \\ 
z^{\left( j\right) \delta }\underset{\left( ji\right) }{P}\underset{}{%
_{\delta }}\overset{}{^{\alpha }}\underset{}{_{\beta }}\underset{}{_{\gamma }%
}=\underset{\left( jj\right) }{P}\overset{}{_{\beta \gamma }^{\alpha }}, & 
& z^{\left( 1\right) \delta }\underset{\left( 2i\right) }{Q}\underset{}{%
_{\delta }}\overset{}{^{\alpha }}\underset{}{_{\beta }}\underset{}{_{\gamma }%
}=\underset{\left( i2\right) }{C}\overset{}{_{\beta \gamma }^{\alpha }}, & 
& z^{\left( 2\right) \delta }\underset{\left( 2i\right) }{Q}\underset{}{%
_{\delta }}\overset{}{^{\alpha }}\underset{}{_{\beta }}\underset{}{_{\gamma }%
}=\underset{\left( 22\right) }{\overset{i}{Q}}\overset{}{_{\beta \gamma
}^{\alpha }}, \\ 
&  &  &  &  \\ 
z^{\left( j\right) \delta }\underset{\left( ji\right) }{S}\underset{}{%
_{\delta }}\overset{}{^{\alpha }}\underset{}{_{\beta }}\underset{}{_{\gamma }%
}=\underset{\left( i\right) }{\overset{j}{S}}\overset{}{_{\beta \gamma
}^{\alpha }}, &  & z^{\left( 1\right) \delta }\underset{\left( 2i\right) }{S}%
\underset{}{_{\delta }}\overset{}{^{\alpha }}\underset{}{_{\beta }}\underset{%
}{_{\gamma }}=0, &  & z^{\left( 2\right) \delta }\underset{\left( 1i\right) }%
{S}\underset{}{_{\delta }}\overset{}{^{\alpha }}\underset{}{_{\beta }}%
\underset{}{_{\gamma }}=\underset{\left( 12\right) }{R}\overset{}{_{\beta
\gamma }^{\alpha }}.%
\end{array}
\tag{5.4}  \label{5.6}
\end{equation}%
$\left( i=0,1,2;j=1,2\right) $

\begin{acknowledgement}
This work was supported by Contract with Sinoptix No. 8441/2009.
\end{acknowledgement}

Alexandru \textsf{OANA}

University Transilvania of Bra\c{s}ov, Department of Mathematics and
Informatics, Blvd. Iuliu Maniu, no. 50, Bra\c{s}ov 500091, Romania.

\textit{E-mails:} alexandru.oana@unitbv.ro

\end{document}